\documentclass[preprint,12pt]{elsarticle}
\usepackage{amssymb}
\usepackage[english]{babel}

\usepackage{amsfonts}
\usepackage{amsmath}
\usepackage{graphicx}
\usepackage[colorlinks=true, allcolors=blue]{hyperref}
\usepackage{natbib}

\newtheorem{prop}{Proposition}
\newtheorem{thm}{Theorem}
\newtheorem{rem}{Remark}

\begin{document}

\begin{frontmatter}
\title{Periodic Solutions in a Simple Delay Differential Equation}

\author[inst1]{Anatoli Ivanov}

\affiliation[inst1]{organization={Department of Mathematics},
            addressline={Pennsylvania State University}, 
            country={USA}}

\author[inst2]{Sergiy Shelyag}

\affiliation[inst2]{organization={College of Science and Engineering},
            addressline={Flinders University}, 
            country={Australia}}





\begin{abstract}
Simple form scalar differential equation with delay and nonlinear negative periodic feedback is considered. 
The existence of several types of slowly oscillating periodic solutions is shown with the same and double periods of the feedback 
coefficient. The periodic solutions are built explicitly in the case of piecewise constant nonlinearities 
involved. The periodic dynamics are shown to persist under small perturbations of the equation which make it smooth. 
The theoretical results are verified by extensive numerical simulations.
\end{abstract}


\end{frontmatter}

\section{Introduction}\label{Intro}

The simple form delay differential equation (DDE)
\begin{equation}\label{DDE}
    x^\prime(t)=-\mu x(t)+f(x(t-\tau))
\end{equation}
represents a broad variety of complex dynamical behaviors that can happen in delay differential models. 
It is also one of the most studied equations in the field exhibiting the range of diverse dynamics from global stability
of equilibria to instability and existence of periodic solutions and chaotic behaviours (in contrast to 
analogous ordinary differential equations). Basic and fundamental facts about the DDE (\ref{DDE}) 
can be found in many sources, see e.g. monographs \cite{Ern09,DieSvGSVLWal95,HalSVL93,Smi11}, among others.

Equation ({\ref{DDE}}) is also extensively used as a mathematical model of a variety of real world phenomena, most
notably in biological applications; see e.g. the monographs \cite{Ern09,GlaMac88,Kua93,Smi11} and further references
therein. Among the well known and most studied models of type (\ref{DDE}) are the Mackey-Glass physiological equations
\cite{GlaMac88,MacGla77}, Wazewska-Lasota blood cell model \cite{WazLas76}, Nicholson's
blowfly equations \cite{GurBlyNis80,PerMalCou78}, and several others \cite{Ern09,Kua93,Smi11}.

Some of the more accurate and adequate mathematical models are explicitly time-dependent, which take into account
certain intrinsic periodicity factors, such as e.g. circadian rhythms and seasonal changes. Examples of such 
differential delay models can be found in publications \cite{AmsIde13,FriWu92,GlaMac88,Kua93,LiKua01,YuLi19}. 
Therefore, a natural generalization of equation (\ref{DDE}) is the following DDE 
\begin{equation}\label{DDE_P}
    x^\prime(t)=-\mu(t) x(t)+a(t)f(x(t-\tau)),
\end{equation}
where $\mu(t)$ and $a(t)$ are $T$-periodic functions, and the nonlinearity $f(x)$ satisfies appropriate feedback conditions. 
An important theoretical and applied question then arises whether the DDE (\ref{DDE_P}) admits 
periodic solutions, which are induced by the periodicity in the coefficients $\mu$ and $a$. Therefore, the mathematical 
problem to address is to derive conditions under which equation (\ref{DDE_P}) possesses a nontrivial  $T$-periodic 
solution. This research work is our initial attempt to answer this question for some specific partial cases of the equation
under consideration.

{  This paper deals with the problem of existence of periodic solutions for the simplest possible scalar differential delay equation of the type (\ref{DDE_P}) when} $\mu(t)\equiv0$,
\begin{equation}\label{DDE_SP}
    x^\prime(t)=a(t)f(x(t-\tau)),
\end{equation}
and where $f(x)$ is a continuous function satisfying the negative feedback assumption $x\cdot f(x)<0,\;\forall x\ne0,$ 
the coefficient $a(t)\ge0$ is a continuous periodic function with the period $T$ satisfying $T>\tau,$ and $\tau>0$ is a delay.
The principal problem we are addressing in this work is to derive conditions on the parameters and nonlinearities of
equation (\ref{DDE_SP}) which would yield the existence of periodic solutions with the same period $T$ as the coefficient
$a(t).$ Though there are some results on the existence of periodic solutions for the DDEs with periodic coefficients 
none of them is applicable to the simplest form equation (\ref{DDE_SP}), {  to the best of our assessment and knowledge}. In many cases that are published in the research literature the generalization from the autonomous case to the period one leads to the loss of the steady positive equilibrium. In some others the  equilibrium state persists, however, the sufficient conditions for the existence of periodic solutions cannot distinguish between the nontrivial periodic solutions and the steady state. {  See more details on such relevant results and further references in e.g. \cite{AmsIde13,Far17,FriWu92,LiKua01,YuLi19}.}

{  We are stating and  dealing with the problem of existence of nontrivial periodic solutions in the simple form DDE (\ref{DDE_SP}). Due to the negative feedback assumption on $f$ this equation always admits the trivial constant solution $x(t)\equiv0$. The mathematical problem is to find periodic solutions which are distinct from the trivial one. This is a difficult problem which is largely not addressed on the general theoretical level. One way to pave a path for a systematic analysis in this direction is to consider model equations of the form given by equation (\ref{DDE_SP}) with simple $f$ and $a$.}

We first construct the periodic solutions in an explicit form by using piecewise constant functions (both the nonlinearity $f$ 
and the coefficient $a$). Then we modify both functions to be continuous and ``close'' to the piecewise constant
ones, by ``smoothing'' them in a small neighborhood of the discontinuity set (they can also be
made of the $C^\infty$ type). The periodic dynamics and their stability are shown to persist under such small perturbations.



The paper is structured as follows. First section ``Preliminaries'' contains basics of the DDE (\ref{DDE_SP}) 
and preliminary results  from our 
recent conference papers \cite{IvaShe23,IvaShe23-2} necessary for the exposition in subsequent sections. Section \ref{PerSols} deals 
with periodic solutions of equation (\ref{DDE_SP}) explicitly constructed for the initially piecewise constant nonlinearity $f$ and 
the coefficient $a$.  The dynamics  on a set of slowly oscillating solutions are reduced to that for induced one-dimensional maps, 
thus proving both the existence of periodic solutions and the types of their stability. In Section 3 it is shown that the dynamics of
the interval maps and the periodic solutions from section 2 persist when the defining functions $f$ and $a$ of equation 
(\ref{DDE_SP}) are replaced by close to them continuous or smoothed functions. The arguments of small regular perturbations of 
dynamics by DDE (\ref{DDE_SP}) are used, which are behind such ``small modifications''. The theoretical derivations and exact 
analytical calculations are verified by extensive numerical simulations.

\section{Preliminaries}\label{Prelim}


We focus in this work on the particular case $\mu=0$ and $\tau=1$ when equation (\ref{DDE}) becomes 
\begin{equation}\label{DDE-M}
    x^\prime(t)=a(t)f(x(t-1)).
\end{equation}
Note that the case of general delay $\tau>0$ is easily reduced to the standard consideration $\tau=1$ by the time scaling $t=\tau s$.

Given an initial function $\varphi(s)\in C([-1,0],\mathbb{R}):= \mathbb{X}$ in the standard phase space $\mathbb{X}$ the 
corresponding solution $x=x(t,\varphi)$ to equation (\ref{DDE-M}) is easily found for all $t\ge0$ by consecutive forward 
integration (the step method). Let $S^{t}(\varphi)$ be the shift in $\mathbb{X}$ by time $t$ along the solutions of delay 
differential equation (\ref{DDE-M}), that is $S^{t}(\varphi)=x(t+s,\varphi), s\in[-1,0]$. It is a  straightforward observation
that a fixed point $\varphi_0$ of the shift by period $T$, $S^{T}(\varphi_0)=\varphi_0$, gives rise to a periodic solution 
$x(t,\varphi_0)$ of the DDE (\ref{DDE-M}). Under such an interpretation, the zero solution $x(t)\equiv0$ is a result of the trivial 
fixed point $\varphi_0(s)\equiv0, s\in[-1,0]$, which always exists for DDE (\ref{DDE-M}), due to the negative feedback assumption. 
The mathematical problem is to derive conditions for the existence of non-trivial fixed points of the shift operator $S^{T}$, 
which is a principal objective of this work.

A solution is called oscillating if it has an infinite sequence  of zeros $\{t_n\}\to\infty, 
x(t_n)=0, n\in\mathbb{N}$, and $x(t)$ is not an identical zero eventually (on any interval of the form $[t_0,\infty$)). 
An oscillating solution is called slowly oscillating, if the sequence of zeros is such that the distance between consecutive zeros 
is greater than the delay, i.e. $t_{n+1}-t_n>1, \forall n\in\mathbb{N}.$

Sufficient conditions for the oscillation of all solutions to equation (\ref{DDE-M}) are well known, see e.g. 
\cite{GyoLad91} and further references therein. With the negative feedback assumption $x f(x)<0, x\ne0,$ and 
the non-negativity of the coefficient $a(t)\ge0, a(t)\not\equiv0,$ they are that either $f^\prime(0)$ is sufficiently
large, $|f^\prime(0)|\ge f_1>0,$ and the coefficient $a$ is separated from zero, $a(t)\ge a_1>0$, or $f^\prime(0)\ne0$
and $a_1>0$ is large enough (in both cases the product $|f_1\cdot a_1|$ must be sufficiently large). These conditions are
satisfied for our considerations to follow in next sections, where the functions $f$ and $a$ are smoothed continuous ones derived 
from the initial piecewise constant functions (the piecewise constant case can be viewed as a limit of the continuous one). 

We deal in this paper with the slowly oscillating solutions to equation (\ref{DDE-M}). The slowly oscillating solutions 
are obtained when the initial functions belong to one of the two sign definite cones.
Introduce the following  two standard subsets of initial functions
\begin{align*}
&K_{-}:=\{\varphi\in\mathbb X\;\vert\; \varphi(s)<0\;~\forall s\in[-1,0] \}\text{ and} \\
&K_{+}:=\{\varphi\in\mathbb X\;\vert\; \varphi(s)>0\;~\forall s\in[-1,0]\}.
\end{align*}

Assuming that all solutions to equation (\ref{DDE-M}) oscillate,
it is a straightforward calculation to verify that for arbitrary $\varphi\in K_{-}$ the corresponding
solution $x(t;\varphi)$ has an increasing sequence of zeros $0<t_1<t_2<t_3< \cdots $ such that 
$t_{k+1}-t_{k}>\tau=1, k\in\mathbb N,$ and 
\begin{equation}\label{t_n}
    x(t)>0\;\; \forall t\in(t_{2k-1}, t_{2k})\quad\text{and}\quad x(t)<0\;\; \forall t\in(t_{2k}, t_{2k+1}).
\end{equation}
In the case when $\varphi\in K_{+}, \varphi(0)=h>0,$ the above inequalities on the consecutive intervals $(t_{2k-1},t_{2k})$ and 
$(t_{2k-1},t_{2k})$ are reversed into opposite.
Therefore, the solution $x(t;\varphi)$ is slowly oscillating. The consecutive intervals $(t_j,t_{j+1}), j\in\mathbb{N}$, are called semicycles.

Based on the well-known facts \cite{GyoLad91} and the reasoning above we have the following statement.
\begin{prop}\label{Prop0}
    Suppose that the nonlinearity $f\in C(\mathbb{R,\mathbb{R}})$ satisfies the negative feedback assumption 
    $x\cdot f(x)<0, \forall x\ne0,$ and $f^\prime(0)<0$. Let the coefficient $0< a(t)\in C([0,\infty),\mathbb{R})$ be $T$-
    periodic and such that $\min\{a(t), t\in[0,T]\}\ge m_0>0$. Assume that $|f^\prime(0)|\cdot m_0>{1}/{e}$. Then for every initial
    function $\varphi\in\mathbb{K_{-}}$ with $\varphi(0)=h<0$ the corresponding solution $x(t,h), t\ge0,$ is slowly 
    oscillating with the sequence of simple zeros $\{t_n\}$ satisfying the inequalities (\ref{t_n}) on the consecutive 
    semicycles. Likewise, for every initial function $\psi\in K_+$ with $\psi(0)=h>0$ the corresponding solution 
    $x(t,\psi), t\ge0,$ is slowly oscillating satisfying the inequalities opposite to those in (\ref{t_n}).
\end{prop}
Additional related theoretical basics on delay differential equations can be found in e.g. 
monographs \cite{DieSvGSVLWal95,HalSVL93}.

\section{Periodic Solutions}\label{PerSols}

In this section the periodic solutions of interest are explicitly constructed based on model piecewise constant nonlinearity $f$
and the periodic coefficient $a$. Let $f(x)=f_0(x)=-\text{sign}(x)$ and define $a=a(t,p_1,p_2,a_1,a_2):=A_0(t)$ as follows
\begin{equation}\label{A0}
a(t) = A_0(t) = \begin{cases}
a_1>0,\;\text{if}\; 0\le t < p_1 \\
a_2>0,\;\text{if}\; p_1\le t < p_1+p_2\\
\text{periodic extension outside}\; [0,p_1+p_2)\;\text{for all}\; t\in\mathbb{R}.
\end{cases}
\end{equation}
Let an initial function $\varphi(s)\in\mathbb{X}$ be such that either $\varphi(s)<0$ or $\varphi(s)>0, \forall s\in[-1,0],$ and 
$\varphi(0)=h\ne 0$. 
Clearly, the corresponding solution $x=x(t,\varphi)$ depends on the value $h$ only, and does not depend of the 
particular values of $\varphi(s)\ne0$ for $s\in[-1,0)$. 
In addition, it is uniquely determined by the parameters $p_1,p_2,a_1,a_2$ of the function $a(t)$.
The solution $x$ is slowly oscillating and piecewise linear (affine) for $t\ge0$.

Like in the case of continuous $f$ and $a$ the corresponding solutions to the DDE (\ref{DDE-M}) with $f=f_0$ and
$a=a_0$ are slowly oscillating with the consecutive zeros $\{t_n\}$ satisfying the inequalities (\ref{t_n}). In this case we do not 
need any additional assumptions on $f_0$ or $a_0$ as in Proposition \ref{Prop0}, as the lower rates of growth or decay of any such 
solution are bounded away from zero (the rates of growth/decay are also bounded from above, due to the boundedness of $f$ 
and $a$). {  However, since the constructed solutions are slowly oscillating, we assume that necessarily the assumption $p_1+p_2>1$ 
must be in place.}

In this section we follow and expand on ideas and explicit constructions from our recent work \cite{IvaShe23,IvaShe23-2}. 
In particular, functions  $f$ and $a$ are of  the same form as in the above papers.
 
\subsection{Periodic solutions with coefficient's period (Type I)}\label{PS-1T}

In this subsection we explicitly construct piecewise affine periodic solutions to equation (\ref{DDE-M}) which are of the  same 
period $T>1$ as the periodic coefficient $a(t)$ is. Some of the exposition elements are close and based on results in our recent 
conference paper \cite{IvaShe23-2}.

\begin{figure*}
    \centering
    \includegraphics[width=0.9\textwidth]{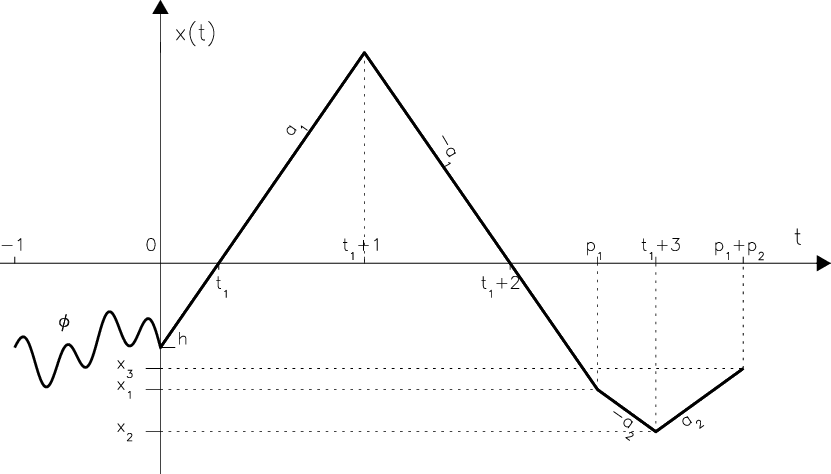}
    \caption{Slowly oscillating piecewise affine solution (Type I) }
    \label{fig:Fig1}
\end{figure*}

We presume that the desired periodic solutions have the shape as shown in Figure \ref{fig:Fig1}. That is, the initial value 
$\varphi(0)=h<0$ is determined by an initial function $\varphi(s)\in K_-$. The corresponding solution $x=x(t,\varphi), t\ge0,$  has
exactly two zeros $\{t_1,t_2\}$ on the initial period $[0,T], T=p_1+p_2,$ with $t_1<t_2=t_1+2<T$, for particular assumptions on the 
parameters $p_1,p_2, a_1,a_2$. Such assumptions are verified to be valid by the analytical calculations that follow, as well as they
are confirmed by corresponding numerical simulations. 
It is straightforward to find desired values of the solution $x=x(t,h)$ for any $t\in [0,T]$. 
We have the following calculations:
$$
t_1=-\frac{h}{a_1},\quad x_1=x(p_1)=-h+2a_1-a_1p_1,
$$
$$
x_2=x(t_1+3)=(\frac{a_2}{a_1}-1)h+a_1(2-p_1)+a_2(p_1-3), 
$$
$$
x_3=x(T)=x(p_1+p_2)= (2\frac{a_2}{a_1}-1)h+a_1(2-p_1)+a_2[(2p_1+p_2)-6]:=G(h),
$$
where the coefficients of the  affine function $G(h):=mh-b$ are given by
\begin{equation}\label{m-b}
 m=2\frac{a_2}{a_1}-1,\quad\text{and}\quad   b=a_1(p_1-2)+a_2[6-(2p_1+p_2)].
\end{equation}

It is elementary to deduce that the range for the slope $m$ is the open interval $(-1,\infty)$ with $|m|<1$ when $0<a_2<a_1$ holds, and $m>1$ when $0<a_1<a_2$ is valid.

A fixed point $h_*<0$ of the map $G(h)=mh-b$ gives rise to a slowly oscillating periodic solution $x=x(t,h_*)$ of 
equation (\ref{DDE-M}). Moreover, by the affine type of $G$ and the construction, such periodic solution is asymptotically 
stable if $|m|<1$,  which is equivalent to  $0<a_2<a_1$. The periodic solution is unstable when $m>1$ , which is equivalent
to $0<a_1<a_2$.
The unique fixed point $x=h_*$ is easily found as 
\begin{equation}\label{h}
    h_*=b/(m-1)=\frac{a_1[a_1(p_1-2)+a_2[6-(2p_1+p_2)]]}{2(a_2-a_1)}.
\end{equation}
Therefore, we arrive at the following statement.
\begin{thm}\label{Thm1}
Let the parameters $a_1,a_2,p_1,p_2$ be given with the values $m$ and $b$ defined by (\ref{m-b}). 
Then 
\begin{itemize}
    \item[(a)] 
In case when  $ |m|<1\,~(0<a_2<a_1)$ and $b>0$ are valid, the DDE (\ref{DDE-M}) has an 
asymptotically stable slowly oscillating $T$-periodic solution. 
     \item[(b)]
In the complimentary case when $m>1\,~(0<a_2<a_1)$ and $b<0$ are 
valid the equation has an unstable slowly oscillating $T$-periodic solution. 
\end{itemize}
The periodic solutions are generated by the initial function $\varphi(s)
\equiv h_*,~ s\in[-1,0],$ where $h_*$ is given by (\ref{h}).
\end{thm}
\begin{rem}
\begin{itemize}
    \item[(a)] 
    The local stability case $|m|<1$ is equivalent to the relationship $a_1>a_2\ge0$. At the same time, in order for the
    fixed point to be negative, one must assume that $b=a_1(p_1-2)+a_2[6-(2p_1+p_2)]>0$. The inequality $p_1>2$ is valid by the 
    construction. For arbitrary fixed $p_1>2, p_2>0, a_2>0$ there exists $a_1^0$ such that the inequalities $|m|<1$ and $b>0$ are 
    satisfied for all $a_1\ge a_1^0.$ Thus, the required assumptions are met for all sufficiently large $a_1.$ This case is treated
    in more detail in our recent paper \cite{IvaShe23-2}.
    \item[(b)] 
    The instability case $m>1$ is equivalent to $0<a_1<a_2$ (note that the case $m<-1$ is not possible for any choices of
    positive $a_1$ and $a_2$). In order to have the fixed point $h_*$ negative, one must assume that $b<0$ is satisfied, that is 
    $a_1(p_1-2)+a_2[6-(2p_1+p_2)]<0$. Since $p_1>2$ the assumption $p_2\ge2$ would imply that $6-(2p_1+p_2)<0$. Therefore, given 
    fixed values of $a_1>0, p_1>2, p_2\ge2$ there exists $a_2^0$ such that $b<0$ for all $a_2\ge a_2^0.$ Thus, the required 
    assumptions are met for all sufficiently large $a_2.$
    \item[(c)] 
    It is clear that due to the continuous dependence of $m$ and $b$ on the parameters $a_1, a_2, p_1, p_2$ the stable or unstable periodic solutions of Theorem \ref{Thm1} exist on an open box in $\mathbb{R}^4$, each side of which is a corresponding open interval about a specific value of each of the four parameters corresponding to a fixed point $h_*$. {  Those boxes can be multiple, as our numerical insight indicates, and they  do not overlap for the  stable and unstable solutions due to the respective opposite relationship between $a_1$ and $a_2$.}
\end{itemize}
\end{rem}
\begin{table}
    \begin{tabular}{|c|c|c|c|c|c|}
\hline    
$~~~~~a_1~~~~~$ & $~~~~~a_2~~~~~$ & $~~~~~p_1~~~~~$ & $~~~~~p_2~~~~~$ & $~~~~~h_*~~~~~$ & $~~~~~T~~~~~$ \\
\hline
1 &	0.25	& 2.5	& 1.5	& -0.25	& 4 \\
2	& 0.5	& 2.5	& 2	& -1/3	& 4.5 \\
2	& 0.25	& 2.5 &	1 &	-4/7	 & 3.5 \\
1	& 0.5	& 3	& 1 &	-0.5	& 4 \\
2	& 1	& 3	& 1.5 &	-0.5 &	4.5 \\
2.5	& 0.5	& 3	& 4 & -0.31 & 	7 \\
3	& 0.5	& 3	& 4.5	& -0.45	& 7.5 \\
4   & 1     &3.5 & 2    & -2  & 5.5 \\
5	& 0.5	& 3	& 3	& -1.94 &	6 \\
5	& 1	& 3	& 2	& -1.88	& 5 \\
{ $\sqrt{10}$} & {  $1/\sqrt{5}$} & {  $\pi$} & {  $e+1$} & {  -1.0602} & {  $\pi+e+1$} \\
\hline
    \end{tabular}
    \caption{Examples of parametric values for which stable $T$-periodic solutions of equation~(\ref{DDE-M}) with coefficient $a$ 
    defined by (\ref{A0}) have been demonstrated numerically.}
    \label{tab:param_table}
\end{table}
The parametric range for which Theorem 1 guarantees the existence of slowly oscillating periodic solutions, either stable or unstable, is large.

Note that due to the symmetry property of the nonlinearity $f$ (oddness) and the procedure of construction of the periodic solutions
$x_{h_*}$, under the assumption of Theorem 1 
there also exists the symmetric to $x_{h_*}$ periodic solution generated by
the initial function $\phi(s)\equiv -h_*>0$. The two periodic solutions are related by $x_{-h_*}(t)\equiv-x_{h_*}(t)$, 
thus, they are symmetric to each other.

A large sample of parametric values has been derived for which numerical periodic solutions have been obtained and confirmed to be 
of the described form.
Table~\ref{tab:param_table} is a sample selection of parametric values resulting in asymptotically stable slowly oscillating periodic solutions given by Theorem 1. The values $h_*<0$
for these periodic solutions are easily found from the explicit formula (\ref{h}). They are also verified to exist numerically.

The next Table 2 is a sample selection of unstable slowly oscillating periodic solutions implied by Theorem 1. 
The values $h_*$ for these periodic solutions are found from formula (\ref{h}) as well. The unstable periodic solutions  cannot be
verified numerically as respective initial functions $\varphi\equiv h_*<0$ actually generate approximate close nearby solutions 
(not the exact ones). 
Those solutions, when calculated in forward time, are repelled for large times by the actual unstable ones and are attracted by stable
periodic solutions of double period $2T$. The latter are studied analytically and numerically and described in detail in the next 
subsection.

\begin{table}
    \begin{tabular}{|c|c|c|c|c|c|}
\hline    
$~~~~~a_1~~~~~$ & $~~~~~a_2~~~~~$ & $~~~~~p_1~~~~~$ & $~~~~~p_2~~~~~$ & $~~~~~h_*~~~~~$ & $~~~~~T~~~~~$ \\
\hline
 0.5	& 5	& 5	& 1	& -1.31	& 6 \\
1	&   5	& 4	& 1	& -1.625 & 5 \\
1	&   7 	& 4 & 1 & -1.5	 & 5 \\
1	&   7	& 3	& 2 & -1.0833 & 5 \\
1.5	&  7 	& 2.5 &	2.5 & -1.3295 &	5 \\
1.5	&  7	& 3.5	& 1.5 & -2.0795 & 	5 \\
1.5	&  7 	& 4	& 3	& -4.3636	& 7 \\
2	&  7 	& 3	& 2	& -2.4 &	5 \\
2	&  7 	& 3	& 3	& -4.6 &	6 \\
2	&  7	& 3.5	& 2.5	& -4.3	& 6 \\
2.5 &  7    & 2.5 &  2.5 & -1.2614 & 5 \\
2.5 &  7    & 3   & 2  &  -1.5682  & 5 \\
{ $1/\sqrt{5}$} & { $3\pi/2$ }& { $2e$ }& { $\pi/3$ }& {  -1.382} & { $2e+\pi/3$}\\
\hline
    \end{tabular}
    \caption{Examples of parametric values for which unstable periodic solutions exist to equation~(\ref{DDE-M}) with the coefficient 
    defined by (\ref{A0}). The values $h_*$ are given by (\ref{h}) with $b<0$ and $m>1$.}
    \label{tab:unstable_table}
\end{table}

\subsection{Periodic solutions with double coefficient's period (Type II)}\label{PS-2T}

In this subsection we explicitly construct piecewise affine periodic solutions to equation (\ref{DDE-M}) which are of 
double period $2T$ of the $T$-periodic coefficient $a(t)$. Some of the exposition is close and based on the results of
our recent conference paper \cite{IvaShe23}.

\begin{figure}
    \centering
    \includegraphics[width=1\textwidth]{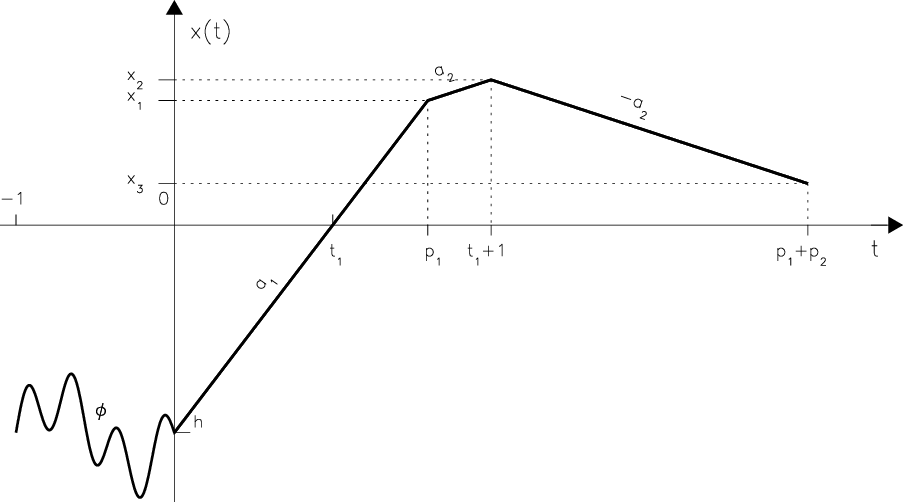}
    \caption{Slowly oscillating piecewise affine solution (Type II) on initial periodic interval}
    \label{fig:Fig2}
\end{figure}

We presume that the desired periodic solutions have the initial shape as shown in Figure \ref{fig:Fig2}. That is, the 
initial value $\varphi(0)=h<0$ is determined by an initial function $\varphi(s)\in K_-$. The corresponding solution 
$x=x(t,\varphi), t\ge0,$  has exactly one zero $t_1$ on the initial periodic interval $[0,T], T=p_1+p_2,$ with $t_1<p_1$
and $x(t)<0 \; \forall t\in[0,t_1)$ and $x(t)>0 \;\forall t\in(t_1,T]$. Such assumptions are verified to be valid by 
the analytical calculations that follow, as well as are confirmed by associated numerical simulations. 
It is straightforward to calculate the desired values of $t_1, x_1, x_2, x_3$ as follows:
$$
t_1=-\frac{h}{a_1},\quad x_1=x(p_1)=h+a_1p_1,\quad x_2=x(t_1+1)=(1-\frac{a_2}{a_1})h +p_1(a_1-a_2)+a_2,
$$
$$
x_3=x(p_1+p_2)=(1-2\frac{a_2}{a_1})h+a_1p_1+a_2(2-2p_1-p_2):=F_1(h)=kh + d,
$$
where
\begin{equation}\label{m-b2}
    k=1-2\frac{a_2}{a_1}\quad\text{and}\quad d=a_1p_1+a_2(2-2p_1-p_2)
\end{equation}
are uniquely defined by the parameters $p_1,p_2,a_1,a_2$.\newline
Likewise, when $\psi\in K_{+}$ and $\psi(0)=h>0$ analogous calculations yield 
$$
x_3=x(p_1+p_2,\psi)=\left(1-2\frac{a_2}{a_1}\right)\,h-a_1p_1-a_2(2-2p_1-p_2):=F_2(h)=kh-d,
$$
with the same $k$ and $d$ as in the expression for $F_1$ given by formula (\ref{m-b2}).
It is easy to see that for arbitrary positive $a_1, a_2$ the respective values of $k$ can only belong to the range 
$(-\infty, 1)$ with $|k|<1$ if $0<a_2<a_1$ and $k<-1$ if $0<a_1<a_2.$

We are interested in parametric values $p_1,p_2,a_1,a_2$ and initial values $h<0$ such that $h_1=F_1(h)>0$ and $h_2=F_2(h_1)<0$.
This will be the case when $d>0$, at least in some vicinity $|h|<\delta_1, h\ne0,$ of the discontinuity point $h=0$. 
The through map $F_0=F_2\circ F_1$ is given by $F_0(h)=k^2 h+(k-1)d$. Its only fixed point is found as $h_*=-d/(k+1)$, and it is 
negative when $d>0$ and $|k|<1.$ The fixed point $h_*$ generates an asymptotically stable periodic solution of period $2T$, when the
initial function for equation (\ref{DDE-M}) is chosen as $\varphi(s)\equiv h_*.$ The periodic solution can also be described in terms of
a cycle of period $2$ of an appropriate one-dimensional map $F$ defined by the maps $F_1$ and $F_2$ introduced above.
Namely, define the piecewise affine map $F$ for $h\in\mathbb{R}\setminus\{0\}$ as follows:
\begin{equation*}
    F(h)=\begin{cases} & F_1(h)=kh+d,\quad\text{if}\quad h<0\\
                       & F_2(h)=kh-d,\quad\text{if}\quad h>0.
    \end{cases}
\end{equation*}
The map $F$ is discontinuous at $h=0$ and is symmetric (odd) with respect to $h=0$, $F(-h)=-F(h),\; \forall h\ne0$. 
Assuming $d>0$ one can see that the negative feedback condition is satisfied at least locally in a vicinity of $h=0$.
The global dynamics of map $F$ is very simple in this case, as described by the following statement:
\begin{prop}\label{Prop2}
    Suppose that $d>0$. Then
    \begin{itemize}
        \item[(i)] If $0<a_2<a_1$ is satisfied, the map $F$ has a unique globally attracting cycle of period $2$, given by 
        $\{h_*,-h_*\}$, where $h_*=-{d}/(k+1)$. The $2$-cycle attracts all the initial values
        $h\in\mathbb{R}\setminus\{0\}$ in the case when $a_1>2a_2$ holds, and all the initial values from the interval 
        $h\in(-d/k,d/k)$ in the case when $2a_2>a_1>a_2$ holds;
        \item[(ii)] If $0<a_1<a_2$ is satisfied, the map $F$ has no cycles of period 2; moreover, any iterative sequence 
        $F^n(h), h\ne0,$ diverges, $\lim_{n\to\infty} |F^{n}(h)|=\infty$.
    \end{itemize}
\end{prop}
Therefore, by translating the properties of the interval map $F$ to the slowly oscillating solutions of DDE
(\ref{DDE-M}), we have the following statement
\begin{thm}\label{Thm2}
    Let the parameters $a_1,a_2,p_1,p_2$ be given with values $k$ and $d$ defined by (\ref{m-b2}) and $d>0$. 
    Then, in case when  $ |k|<1\, (0<a_2<a_1)$ is valid DDE (\ref{DDE-M}) has an 
    asymptotically stable slowly oscillating $2T$-periodic solution.  The periodic solution is generated by the initial function 
    $\varphi(s)\equiv h_*, s\in[-1,0],$ where $h_*=-{d}/(k+1)$. In the case when $k<-1\,(0<a_2<a_1)$ is valid
    the equation does not have any $2T$-periodic slowly oscillating solutions.
\end{thm}
Table 3 provides sample selections of parametric values for which the stable $2T$-periodic solutions exist. 
Besides using formula (\ref{h}) to calculate respective $h_*$, their existence was also confirmed numerically.
\begin{table}
    \centering
    \begin{tabular}{|c|c|c|c|c|c|}
\hline    
$~~~~~a_1~~~~~$ & $~~~~~a_2~~~~~$ & $~~~~~p_1~~~~~$ & $~~~~~p_2~~~~~$ & $~~~~~h_*~~~~~$ & $~~~~~2T~~~~~$ \\
\hline
4 & 1 & 0.5 & 2.5 & -0.33 & 6 \\
 4	& 1	& 1	& 3.5	& -0.33	& 9 \\
 5	& 1	& 1	& 3.5	& -0.938	& 9 \\
 5 & 1 & 0.5 & 3 & -0.313  & 7 \\
6	& 1	& 1	& 3.5	& -1.5	& 9 \\
6	& 1	& 1	& 4.5	& -0.9	& 11 \\
6 & 1.5 & 1.0 & 3.5 & -0.5 & 9 \\
7 & 1.5 & 1.5 & 3.5 & -2.39 & 10 \\
7 & 2.5 & 2.5 & 2.5 & -2.92 & 10 \\
7 & 2.5 & 2 & 3 & -1.17 & 10 \\
{ $1/\sqrt{5}$} & { $3\pi/2$ }& { $3\pi/2$} & { $e/2$} & { -1.12} & { $3\pi+e$} \\
\hline
    \end{tabular}
    \caption{Examples of parametric values for which stable periodic solutions with double period of coefficient $a(t)$ exist to equation~(\ref{DDE-M}). The values of $h_*$ are given by formula~(\ref{m-b2}).}
\end{table}

\subsection{Coexistence of periodic solutions}\label{CPS}

Any slowly oscillating solution considered above in subsections \ref{PS-1T} and \ref{PS-2T} is uniquely defined by the sign 
definite initial function/values $\varphi\in\mathbb{X}$ such that $\varphi(0)=h\ne0$, and  $\varphi(s)\ne0\;\forall s\in[-1,0],$ 
and by the parameters $p_1,p_2,a_1,a_2$ of the periodic coefficient $a(t)$ (note that $f=f_0(x)=-\text{sign}\,(x)$ is fixed). 
One solves the initial value problem for $t\ge0$:
\begin{equation}\label{IVP-0}
    x^\prime(t)=a_0(t)\,f_0(x(t-1)),\;t\ge0,\; \varphi\in\mathbb{X}, \varphi(0)=h\ne0, \varphi(s)\ne0\; \forall s\in[-1,0].
\end{equation}
Its solution $x_{\varphi}(t)=x_{\varphi}(t,h,p_1,p_2,a_1,a_2)$ is easily found to exist for all $t\ge0$; it is uniquely determined
by the set of the parameters $p_1,p_2,a_1,a_2$ and the initial value of $h$ at $t_0=0$.

Instead of the initial time value $t_0=0$ as in the case above, one can consider a similar initial value problem for any $t_0>0$. 
In particular, in the case of the initial time $t_0=p_1$ one is solving the following initial value problem for $t\ge p_1$:
\begin{equation}\label{IVP-N}
    x^\prime(t)=a_0(t)\,f_0(x(t-1)),\; t\ge p_1,\; \psi(p_1)=h\ne0,~\psi(s)\ne0\;~\forall s\in[p_1-1,p_1].
\end{equation}
Similarly, its solution $x_{\psi}(t)=x_{\psi}(t,h,p_1,p_2,a_1,a_2)$ is uniquely defined by the set of the parameters 
$p_1,p_2,a_1,a_2$ and by the initial value $h$ at $t_0=p_1$.

Now, due to the $T$-periodicity of the coefficient $a(t)$, one can easily see that the solutions of the two initial value problems
(\ref{IVP-0}) and (\ref{IVP-N}) are closely related, with an interchange in the set of the parameters and a time shift. 
Namely, they are just a time shift by $p_1$ with the parameter interchange in the following sense:
\begin{equation}\label{Shift}
    x_{\varphi}(t,h,p_2,p_1,a_2,a_1)=x_{\psi}(t+p_1,h_1,p_1,p_2,a_1,a_2),\quad \forall t\ge0,
\end{equation}
with an appropriate $h_1=h_1(h).$
This simple observation allows to conclude the simultaneous existence of periodic solutions of both types I and II, which are 
described in Theorems \ref{Thm1} and \ref{Thm2}. Since the value of the ratio $a_1/a_2$ with respect to $1$ determines the stability
of those periodic solutions, their stability types must be opposite to each other. We have the following statement
\begin{thm}\label{Thm3}
    There is an open set of parameters $p_1,p_2,a_1,a_2$ such that the DDE (\ref{DDE-M}) has simultaneously an
    unstable slowly oscillating $T$-periodic solution of type I and an asymptotically stable $2T$-periodic solution of type II. 
    Such two periodic solutions can be paired in the following sense: small perturbations of the unstable solution are repelled by
    it and attracted by a respective stable solution, as $t\to\infty.$
\end{thm}

\begin{figure}
    \includegraphics[width=0.9\textwidth]{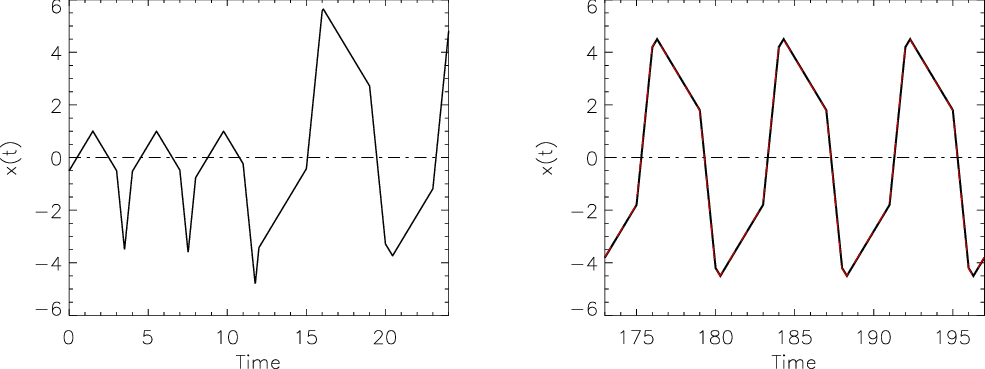}
    \caption{An example of numerical solution with $a_1=1$, $a_2=6$, $p_1=3$, $p_2=1$. 
    Left panel - unstable solution with single period, which is destroyed after two first periods. Right panel - stable solution with double period. Numerical solution of equation~\ref{DDE-M} with $a_1=6$, $a_2=1$, $p_1=1$, $p_2=3$, whose solution has double period of the $T$-periodic coefficient, is overplotted as a red dashed line in the right panel.}
    \label{fig:Fig5}
\end{figure}

Such situation of coexisting two periodic solutions is depicted in Figure \ref{fig:Fig5}. The unstable periodic solution $x_u(t)$ of
period $T=4$ is given by Theorem \ref{Thm1} when $p_1=3,p_2=1,a_1=1,a_2=6$; the respective value of $h$ is $h_*=-0.5$. The stable 
periodic solution $x_s(t)$ of period $2T=8$ is given by Theorem \ref{Thm2} for the parameter values $p_1=1,p_2=3,a_1=6,a_2=1$, 
so that the respective values of $p_1,p_2$ and $a_1,a_2$ are interchanged for these two solutions. In addition, the solution 
$x_s(t)$ is shifted by $p_1$ to the right (the graph (b) in Figure \ref{fig:Fig5} shows the solution $x_s(t+p_1)$ for longer time
intervals).

Tables~\ref{tab:stable_double_table1} and \ref{tab:stable_double_table2} provide a sample selection of such dual cases of the 
coexistence of two types of stable/unstable solutions. Table \ref{tab:stable_double_table2} is a mirror image of Table 
\ref{tab:stable_double_table1} for the stable periodic solutions with double periods. Namely, each row in Table 
\ref{tab:stable_double_table2} represents a stable periodic solution of the corresponding unstable solution in the same row number 
from Table \ref{tab:stable_double_table1}.

\begin{table}
    \begin{tabular}{|c|c|c|c|c|c|}
\hline    
$~~~~~a_1~~~~~$ & $~~~~~a_2~~~~~$ & $~~~~~p_1~~~~~$ & $~~~~~p_2~~~~~$ & $~~~~~h_*~~~~~$ & $~~~~~T~~~~~$ \\
\hline
0.5	& 2.5 & 3	& 0.5	& -0.09	& 3.5 \\
0.5	& 3	& 5	& 1	& -1.35	& 6 \\
0.5	& 5	& 4	& 0.5 &	-0.64 & 4.5 \\
0.5	& 7 & 3	& 2	& -0.52	& 5 \\
1   & 6 & 3 & 1	& -0.5  & 4 \\
1	& 7	& 5	& 1	& -2.67	& 6 \\
2	& 7	& 2	& 3	& -1.4	& 5 \\
3	& 7	& 4	& 1	& -5.62	& 5 \\
{ $\pi/6$} &{  $e$} & { $\pi$} & { 1/2} &{  -0.18} & { $\pi+1/2$} \\
\hline
    \end{tabular}
    \caption{Examples of parametric values for which unstable and stable periodic solutions coexist to DDE~(\ref{DDE-M}) with the coefficient defined by (\ref{A0}). The parametric data is shown for the unstable solutions.}
    \label{tab:stable_double_table1}
\end{table}

\begin{table}
    \begin{tabular}{|c|c|c|c|c|c|}
\hline    
$~~~~~a_1~~~~~$ & $~~~~~a_2~~~~~$ & $~~~~~p_1~~~~~$ & $~~~~~p_2~~~~~$ & $~~~~~h_*~~~~~$ & $~~~~~2T~~~~~$ \\
\hline
2.5	& 0.5 & 0.5	& 3	& -0.156	& 7 \\
3	& 0.5	& 1	& 5	& -0.3	& 12 \\
5	& 0.5	& 0.5	& 4 &	-0.56 & 9 \\
7	& 0.5 & 2	& 3	& -6.19	& 10 \\
6   & 1 & 1 & 3	& -1.8  & 8 \\
7	& 1	& 1	& 5	& -1.17	& 12 \\
7	& 2	& 3	& 2	& -6.3	& 10 \\
7	& 3	& 1	& 4	& -4.375	& 10 \\
{ $e$} & { $\pi/6$} & { 1/2} &  { $\pi$} & { -1.92} & { $2\pi+1$} \\
\hline
    \end{tabular}
    \caption{Examples of parametric values for which stable periodic solutions exist to equation~(\ref{DDE-M}) with the 
    coefficient defined by equation~(\ref{A0}). The parametric data is shown for the stable solutions.}
    \label{tab:stable_double_table2}
\end{table}

\section{Smoothed nonlinearities}\label{Smooth}

In this section, we follow a standard procedure of replacing the discontinuous piecewise constant functions $f_0$ and $A_0$ by close to them 
continuous functions $f_\delta$ and $A_\delta$ for small $\delta>0$.
This procedure has been used in numerous cases; typical examples can be found in papers {  \cite{HalLin86,IvaSha91,Pet80,Wal81b}}.
In a $\delta$-neighborhood of every discontinuity point for both $f_0$ and 
$A_0$ the jump discontinuity is replaced by an affine function by connecting the two constant levels by a line segment.  That is, define $f_\delta(x)$
and $A_\delta(t)$ as follows
\begin{equation}\label{f_d}
    f(x)=f_{\delta}(x)=
\begin{cases}
+1\; \text{if}\; x\le -\delta\\
-1\; \text{if}\; x\ge \delta\\
-({1}/{\delta}) x\; \text{if}\; x\in(-\delta, \delta),
\end{cases}
\end{equation}
and
\begin{equation}\label{A0d}
a(t)=A_{\delta}(t)=\begin{cases}
a_2+\frac{a_1-a_2}{2\delta}(t+\delta)\;\text{if}\; t\in(-\delta,\delta)\\
a_1\;\text{if}\; t\in[\delta,p_1-\delta]\\
a_1+\frac{a_2-a_1}{2\delta}[t-(p_1-\delta)]\;\text{if}\; t\in(p_1-\delta,p_1+\delta)\\
a_2\;\text{if}\; t\in[p_1+\delta,p_1+p_2-\delta]\\
a_2+\frac{a_1-a_2}{2\delta}[t-(p_2-\delta)]\;\text{if}\; t\in(p_1+p_2-\delta,p_1+p_2+\delta)\\
\text{periodic extension on}\; \mathbb{R}\; \text{outside interval}\; [0,p_1+p_2).
\end{cases}
\end{equation}

We shall show next that the values of the piecewise affine slowly oscillating solutions constructed in Section \ref{PerSols} for 
$t\ge0$ do not change on the entire periodic interval $[0,T]$, except in a small $\varepsilon$-vicinity of the corner points 
($\varepsilon$ is a multiple of $\delta$). 

\begin{figure}
    \centering
    \includegraphics[width=0.95\textwidth]{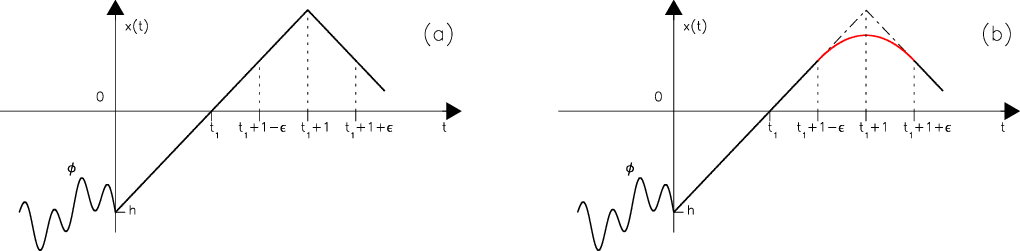}
    \caption{Sketches of (a) continuous solution $x$ with a corner-type discontinuity for $x'$ at $t_1+1$ for functions $A_0(t)$ and $f_0(x)$. (b) $C^1$-smooth solution with parabolic matching on $[t_1+1-\varepsilon,t_1+1+\varepsilon], \varepsilon a_1=\delta,$ for functions $A_\delta(t)$ and $f_\delta(x)$ with small $\delta>0$.}
    \label{fig:Fig3}
\end{figure}
\begin{figure}
    \includegraphics[width=0.95\textwidth]{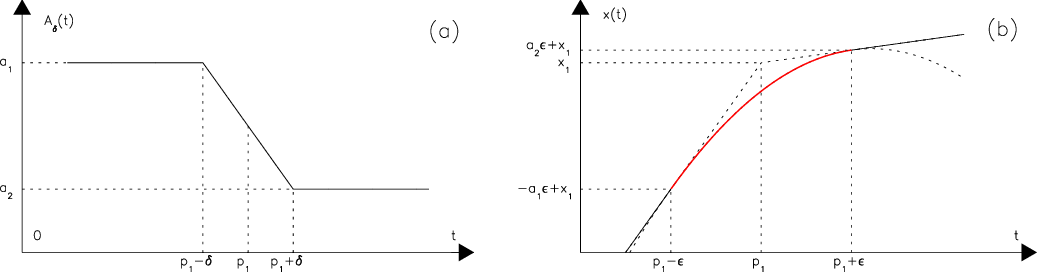}
    \caption{Sketches of (a) continuous coefficient $A_\delta(t)$ in $\delta$-vicinity of $t_1=p_1$; 
    (b) $C^1$-smooth solution with parabolic matching on $[p_1+1-\varepsilon,p_1+1+\varepsilon ]$ for functions $A_\delta(t)$ and $f_\delta(x)$ with small $\delta>0$.}
    \label{fig:Fig4}
\end{figure}

Indeed, the case of a symmetric corner solution, as it passes through a $\delta$-neighborhood of $x=0$ for the nonlinearity 
$f_\delta$, is shown in Figure \ref{fig:Fig3} (which corresponds to e.g. the Type I solution on the first positive semi-cycle 
depicted in Figure \ref{fig:Fig1}). The $\delta$-neighborhood of $x=0$ for $f$ corresponds to the $\varepsilon=\delta/a_1$ 
neighborhood for $t$ at $t=t_1+1$. It is immediately clear that the replacement of the corner type solution on the interval 
$[t_1+1-\varepsilon,t_1+1+\varepsilon]$ by the respective parabolic segment, as a result of the corresponding  integration
now with a linear nonlinearity $f$, results in a symmetric parabola with the vertex at $t_1+1$ and the same equal values 
$x(t_1+1-\varepsilon)$  and $ x(t_1+1-\varepsilon)$ as for the initial piecewise solution with the discontinuous $f_0$. 
The general case of a corner type solution (piecewise affine) is shown in Figure \ref{fig:Fig4}, which corresponds to e.g.
point $(p_1,x_1)$ on the graph of a solution of type II (see Figure \ref{fig:Fig2}). We shall show now that a smooth 
parabolic segment solution connecting points $(p_1-\varepsilon, x((p_1-\varepsilon))$ and 
$(p_1+\varepsilon, x((p_1+\varepsilon))$ is the precise replacement of the corresponding piecewise affine segment solution 
which is obtained for the piecewise constant functions $f_0$ and $A_0$.  The values of the solution $x$ and its derivatives
for the piecewise affine solution are given as 
$$
x(p_1-\varepsilon)=x_1-a_1\varepsilon,\quad x^\prime(p_1-\varepsilon)=a_1,\quad x(p_1+\varepsilon)=x_1+a_2\varepsilon,\quad
x^\prime(p_1+\varepsilon)=a_2,\; 
$$
Looking for a parabolic connection $y=P(t)=A(t-p_1)^2+B(t-p_1)+C$ as a solution of the DDE (\ref{DDE-M}) 
on the interval $[p_1-\varepsilon,p_1+\varepsilon]$, starting at the point $(p_1-\varepsilon,x_1-a_1\varepsilon)$ and 
satisfying the necessary initial-boundary conditions:
$$
x(p_1-\varepsilon)=x_1-a_1\varepsilon,\quad x^\prime(p_1-\varepsilon)=a_1,\quad x^\prime(p_1+\varepsilon)=a_2,\; 
$$
one easily finds that the values of the parameters $A,B,$ and $C$ are:
$$
A=\frac{a_2-a_1}{4\varepsilon},\quad B=\frac{a_1+a_2}{2},\quad C=\frac{a_2-a_1}{4}\,\varepsilon+x_1.
$$
But then, by checking the value of the parabolic segment of the solution $x(t)=P(t)$ at $t=p_1+\varepsilon$, one verifies that
\begin{align*}
P(p_1+\varepsilon)&=A(t-p_1)^2+B(t-p_1)+C\vert_{t=p_1+\varepsilon} \\
            &=\frac{a_2-a_1}{4\varepsilon}\varepsilon^2+
\frac{a_1+a_2}{2}\varepsilon+\frac{a_2-a_1}{4}\varepsilon+x_1=x_1+a_2\varepsilon=x_1+a_2(t-p_1)\vert_{t=p_1+\varepsilon},
\end{align*}
which is the same value as the piecewise affine solution at that point.

Thus, the solutions to equation (\ref{DDE-M}) with continuous $f=f_\delta$ and $a=A_\delta$, are
obtained from the solutions when $f=f_0$ and $a=A_0$ by simple cutting out and replacement of the small corner segments by respective parabolic 
segments, in their $\varepsilon$-vicinity ($\varepsilon=c\cdot\delta$). 
An immediate consequence of this fact is that the periodic solutions and their stability persist under such 
small $\delta$-perturbations of $f_0$ and $a_0$ in DDE (\ref{DDE-M}).

Indeed, the easiest way to see this is to consider the new starting point for $t, t=\varepsilon,$ with $\varphi(\varepsilon)=H\ne0$.
Then the corresponding shift map $S^{T}$ by the period $T$ is dynamically equivalent (conjugate) to the respective map constructed in
the case of $f_0$ and $a_0$. But the latter is equivalent to a composition of respective maps $G,F,F_1.F_2,F_0$ derived in Section 
\ref{PerSols}.
If one would like to stay with the initial value $H:=\varphi(0)\ne0$ in the case of continuous $f_\delta,a_\delta,\delta>0,$ then a 
small adjustment in the considerations should be made. Namely, the new $h:=H\pm c\varepsilon$ should be introduced, where 
$c\varepsilon$ is the difference between the two solutions of the differential delay equation (\ref{DDE-M}) in the two cases when 
$\delta=0$ and small $\delta>0$. Such difference is shown in the graph in Figure \ref{fig:Fig3}, part (b) and in the graph in Figure
\ref{fig:Fig4}, part (b) (in both graphs it is positive). In general, such a difference is independent of the value $x_1$ 
(or equivalent values in other analogous corner situations). Therefore, all the resulting maps $G,F,F_0,F_1,F_2$ derived in Section 
\ref{PerSols} will become finite compositions of the affine maps of the form $h\mapsto\tilde{F},$ where $\tilde{F}
(h)=G_*(h+c_1\varepsilon)+c_2\varepsilon$ and $G_*(h)$ is one of the elementary maps considered in Subsections \ref{PS-1T} and \ref{PS-2T}. 
Therefore, we are in a position to state the following
\begin{thm}\label{Thm4}
    Suppose that the parametric values $a_1, a_2, p_1, p_2$ are such that the assumptions of either one of Theorems \ref{Thm1}
    or \ref{Thm2} are satisfied. Then there exists $\delta_0>0$ such that for every $0\le\delta\le\delta_0$ the delay 
    differential equation (\ref{DDE-M}) with functions $f=f_\delta(x)$ and $a=A_\delta(t)$ has a $C^1$-smooth periodic solution with 
    the same type of stability as in the respective Theorems \ref{Thm1} or \ref{Thm2}. Such a solution converges in uniform metric
    as $\delta\to+0$ to the respective periodic solution with the corner type discontinuity for $x^\prime$ 
    (when $f=f_0, a=A_0$).
\end{thm}
Note that the smoothing of both functions $f_0$ and $a_0$ in a $\delta$-neighborhood of their discontinuity set can be made 
such that the resulting continuous nonlinearities $f_\delta$ and $a_\delta$ are continuously differentiable or even of the 
$C^{\infty}$ class, and Theorem \ref{Thm4} still remains valid. 
Such $C^1$-smoothness is required in some cases
when e.g. one would need to study Floquet multipliers of the periodic solutions, and consider the corresponding linearized
equation along the periodic solutions.

Indeed, it is an elementary  fact that there are many 
solutions for the following approximation problem.  Find a smooth function $F(t)$ defined on an interval $[t_1,t_2]$ 
such that it satisfies the given boundary conditions: $F(t_1)=x_1, F^\prime(t_1)=m_1$ and $F(t_2)=x_2, F^\prime(t_1)=m_2$. 
One of the possible solutions can be suggested as a polynomial of degree 3 or higher.   
In our case of $C^\infty$ smooth connection for $f_0(x)$ on the interval $[-\delta, \delta]$ one can use the following function
$$
f^0(x)=e^{\frac{\delta x}{x-\delta}}-1,\;\text{for}\; x\in[0,\delta],\quad\text{and}\quad f^0(x)=-f^0(x),\; \text{for}\;
x\in[-\delta,0].
$$
It is straightforward to verify that such connecting function is decreasing on the interval $[-\delta,\delta],$ is of class $C^\infty$
there, with all the matching derivatives at $x=0, x=-\delta, x=\delta$. Analogous $C^\infty$-connection can be used for the 
coefficient $a(t)$ at points $0,p_1,T=p_1+p_2$.

\section{Acknowledgements}
The authors thank the mathematical research institute MATRIX in Australia where part of this research was performed. Its final version resulted from collaborative activities of the authors during the workshop ``Delay Differential Equations and Their Applications'' 
(\url{https://www.matrix-inst.org.au/events/delay-differential-equations-and-their-applications/}) held in December 2023. The authors are also grateful for the financial support provided for these research activities by Simons Foundation (USA), Flinders University (Australia), and the Pennsylvania State University (USA). A.I.'s research was also supported in part by the Alexander von Humboldt Stiftung (Germany) during his visit to Justus-Liebig-Universit\"{a}t, Giessen, in June-August 2023.



\end{document}